\documentclass[11pt]{article}

\usepackage{textcomp}
\usepackage{amsmath,amsthm,amssymb,amscd}
\usepackage{cite}

\numberwithin{equation}{section}

\theoremstyle{definition}\newtheorem{definition}{Definition}[section]
\theoremstyle{plain}\newtheorem{theorem}[definition]{Theorem}
\theoremstyle{plain}\newtheorem{proposition}[definition]{Proposition}
\theoremstyle{plain}\newtheorem{lemma}[definition]{Lemma}
\theoremstyle{definition}\newtheorem{assumption}[definition]{Assumption}
\theoremstyle{definition}\newtheorem{example}[definition]{Example}
\theoremstyle{definition}\newtheorem{remark}[definition]{Remark}

\newcommand{\cS}{{\mathcal S}}
\newcommand{\cR}{{\mathcal R}}
\newcommand{\cT}{{\mathcal T}}
\newcommand{\cB}{{\mathcal B}}
\newcommand{\cO}{{\mathcal O}}
\newcommand{\co}{o}
\newcommand{\argmin}{{\mathrm{argmin}}}
\newcommand{\diff}{{\,\mathrm{d}}}
\newcommand{\mybullet}{\mspace{1.5mu}{^{_{_{_\bullet}}}}\mspace{-2.0mu}}

\begin{document}

\title{\bf A new approach to source conditions in regularization with general residual term}

\author{{\sc Jens Geissler}\footnote{Department of Mathematics, Chemnitz University of
Technology, D-09107 Chemnitz, Germany.  E-mail
addresses:$\;$\texttt{jens.geissler\,@\,mathematik.tu-chemnitz.de\,;}
$\;$\texttt{hofmannb\,@\,mathematik.tu-chemnitz.de}.
Research was supported by Deutsche Forschungsgemeinschaft (DFG) under Grant
HO1454/7-2. The work was partly conducted during the Mini Special Semester on Inverse
Problems, May 18 -- July 15, 2009, organized by RICAM (Johann Radon Institute for
Computational and Applied Mathematics, Austrian Academy of Sciences) Linz, Austria.
}
\,\,\,\,\,\,and $\;\;$ {\sc Bernd Hofmann}\footnotemark[1]
 }

\maketitle

\begin{abstract}
This paper addresses Tikhonov like
regularization methods with convex penalty functionals for solving
nonlinear ill-posed operator equations formulated in Banach or,
more general, topological spaces. We present an approach for proving convergence rates which combines advantages of approximate source
conditions and variational inequalities. Precisely, our technique
provides both a wide range of convergence rates and the capability
to handle general and not necessarily convex residual terms as well as nonsmooth operators. Initially formulated for
topological spaces, the approach is extensively discussed
for Banach and Hilbert space situations, showing
that it generalizes some well-known convergence rates results.
\end{abstract}

\section{Introduction}

In recent years because of numerous applications which occurred in
imaging, natural sciences, engineering, and mathematical finance  a
growing interest in different forms of regularization methods for
solving nonlinear ill-posed inverse problems in a Banach space setting could
be observed. This also led to new ideas for proving convergence
rates of such methods in Banach spaces (see,
e.g.,~\cite{Bremen08,BurOsh04,Hein08,Hein09,HeiHof09,Hof09a,HofKalPoeSch07,
KSS09,LorTre08,Neub09,Ramlau08,Resme05,ResSch06,SchGraGroHalLen09,SchLouSch06,Poe08}).
The main problem of handling ill-posed problems in Banach spaces is
the absence of spectral theoretic tools including generalized source
conditions with arbitrary index functions applied to the forward
operator, which were essential for proving results in the Hilbert
space setting.

One way for obtaining convergence rates similar to the well-known
Hilbert space results is the idea of so-called approximate source
conditions, which was originally developed for linear ill-posed
problems in \cite{Hof06} (see also \cite{DHY07,HDK06}) and extended
to nonlinear problems in Banach spaces in \cite{HeiHof09}.
Approximate source conditions, however, rely heavily on the
traditional residual structure being a $p$-th power of the discrepancy
norm. Therefore they are not suited for investigating convergence
rates of variational regularization methods with general residual
terms using appropriate similarity measures. Such progressive
variants of variational regularization were suggested,
comprehensively analyzed, and motivated by means of concrete
examples in \cite{Poe08}.

A second approach, which uses variational inequalities for proving
convergence rates, was first formulated in \cite{HofKalPoeSch07} and
has also been extended to general residual terms in \cite{Poe08}. The
drawback of this second approach in its original form is its
limitation to the standard convergence rate $\cO(\delta)$ for noise
level $\delta>0$ when the reconstruction error is measured by a
Bregman distance.

In this paper, which is mainly based on the thesis \cite{Geis09},
we present an alternative concept that allows both a
wide range of different convergence rates and the use of general
residual terms. Moreover, we  give some new
insight into the interplay of source conditions and variational
inequalities by extending the ideas of \cite{Hof09a} to a more
general setting. Furthermore, we address the question concerning the role
and admissible intervals of an exponent $p>0$ imposed on the
residual term in Tikhonov type regularization (see also
\cite{Hof09b}).

The paper is organized as follows: In Section \ref{sec:notation} we
introduce a Tikhonov type regularization method for the stable
approximate solution of nonlinear ill-posed operator equations in
topological vector spaces with focus on Banach spaces. We formulate
basic assumptions ensuring well-definedness, stability and
convergence of the method and we give a short discussion on
fundamental differences between the classical source conditions and variational
inequalities. In Section \ref{sec:vari} we extend the concept of
variational inequalities introduced in \cite{HofKalPoeSch07}.
Moreover,  we formulate a first convergence rate result in that
section. Based on Section~\ref{sec:vari} in Section~\ref{sec:appvari}, 
which is the main section of this paper, we
present the new approach of approximate variational inequalities for
proving convergence rates of variational regularization methods with
general residual term. In Section~\ref{sec:source} we
restrict our investigations to Banach spaces to clarify the
interplay of approximate source conditions and approximate
variational inequalities. The final Section~\ref{sec:conclude} is devoted to some
concluding remarks where also open questions are formulated.

\section{Problem, notation, and basic assumptions}\label{sec:notation}

Let $F:D(F)\subseteq U\rightarrow V$ be an in general nonlinear operator
possessing the domain $D(F)$ and
mapping between a real topological vector space $U$ and a topological space $V$
with topologies $\tau_U$ and $\tau_V$.
We are going to study operator equations
\begin{equation}\label{eq:fu_v0}
F(u)=v^0
\end{equation}
expressing inverse problems with exact data $v^0\in V$ on the right-hand side.

To ensure mathematical correctness some technical conditions on $U$ and $V$
(in particular the Hausdorff property) are required, but for the important case that $U$
and $V$ are Banach spaces these conditions are always fulfilled.
The topologies $\tau_U$ und $\tau_V$ should be regarded as ``weak'' topologies
because as we will see later in Banach space settings they have to be weaker than the
norm topologies. For this reason we denote convergence with respect to the topology
$\tau_U$ or $\tau_V$ by ``$\rightharpoonup$''.

Instead of the exact right-hand side $v^0$ in (\ref{eq:fu_v0}) only
noisy data $v^\delta$ for some noise level $\delta>0$ are available.
To clarify the meaning of $\delta$ we introduce a non-negative
similarity functional $\cS:V\times V\rightarrow[0,\infty]$, which
not necessarily has to have metric properties, and demand
\begin{equation} \label{eq:simnorm}
\cS(v^\delta,v^0)\leq\delta\qquad\text{and}\qquad\cS(v^0,v^\delta)\leq\delta.
\end{equation}
As approximate solutions of (\ref{eq:fu_v0}) we consider minimizers $u_\alpha^\delta$
over $D(F)$ of the Tikhonov type functional
\begin{equation}\label{eq:tikh}
\cT_\alpha^\delta(u):=\cS(F(u),v^\delta)^p+\alpha\Omega(u)
\end{equation}
with a stabilizing functional $\Omega:U\rightarrow[0,\infty]$, a
regularization parameter $\alpha>0$ and a prescribed exponent
$0<p<\infty$. We set $$D(\Omega):=\{u\in U:\Omega(u)<\infty\} \qquad
\mbox{and} \qquad D:=D(F)\cap D(\Omega).$$

Throughout this paper we make the following assumptions.

\begin{assumption}\label{as:basic}
\mbox{}
\begin{itemize}
\item[(i)]
$F:D(F)\subseteq U\rightarrow V$ ist sequentially continuous with respect to $\tau_U$
and $\tau_V$, i.e. $u_k\rightharpoonup u$ with $u,u_k\in D(F)$ implies
$F(u_k)\rightharpoonup F(u)$.
\item[(ii)]
$D(F)$ is sequentially closed with respect to $\tau_U$, i.e $u_k\rightharpoonup u$ with
$u_k\in D(F)$ and $u\in U$ implies $u\in D(F)$.
\item[(iii)]
There exists a $u\in D$ with $F(u)=v^0$, in particular $D\neq\emptyset$.
\item[(iv)]
The following assertions are fulfilled by $\cS$ (with sequences $(v_k)_{k\in\mathbb{N}}$
and $(\tilde{v}_k)_{k\in\mathbb{N}}$ in $V$ and $v,\tilde{v}\in V$):
  \begin{itemize}
  \item[(a)]
  $\cS(v,\tilde{v})=0$ if and only if $v=\tilde{v}$.
  \item[(b)]
  There exists a value $s\geq 1$ with
  \begin{equation} \label{eq:striangle}
  \cS(v_1,v_2)\leq s\cS(v_1,v_3)+s\cS(v_3,v_2)
  \quad \mbox{for all} \quad v_1,v_2,v_3\in V.
  \end{equation}
  \item[(c)]
  $\cS$ is sequentially lower semi-continuous with respect to $\tau_V$,
  i.e. if $v_k\rightharpoonup v$ and $\tilde{v}_k\rightharpoonup\tilde{v}$ then
  $\cS(v,\tilde{v})\leq\liminf_{k\to\infty}\cS(v_k,\tilde{v}_k)$.
  \item[(d)]
  $\cS(v_k,v)\to 0$ implies $v_k\rightharpoonup v$.
  \item[(e)]
  If $\cS(v_k,v)\to 0$, $\cS(v,v_k)\to 0$, and $\cS(\tilde{v},v)<\infty$ then
  $\cS(\tilde{v},v_k)\to\cS(\tilde{v},v)$.
  \end{itemize}
\item[(v)]
$\Omega$ is convex.
\item[(vi)]
$\Omega$ is sequentially lower semi-continuous with respect
to $\tau_U$, i.e. $u_k\rightharpoonup u$ implies
$\Omega(u)\leq\liminf_{k\to\infty}\Omega(u_k)$.
\item[(vii)]
For each $\alpha>0$ and each $c>0$ the level sets
\begin{equation}\label{eq:levelsets}
M_\alpha(c):=\{u\in D:\cT_\alpha^0(u)\leq c\}
\end{equation}
are sequentially pre-compact with respect to $\tau_U$,
i.e. each sequence $(u_k)_{k\in\mathbb{N}}$ in $M_\alpha(c)$
has a subsequence which converges with respect to $\tau_U$.
\end{itemize}
\end{assumption}

In the sequel for simplicity we will use the terms ``continuous'',
``closed'', and so on instead of ``sequentially continuous'',
``sequentially closed'', and so on if no confusion is to be
expected.
By $U^\ast$ we denote the dual space of $U$, i.e. $U^\ast$ is the set
of all $\tau_U$-continuous linear functionals on $U$.
For $\xi\in U^\ast$ and $u\in U$ we write $\xi(u)$ if we evaluate the functional
$\xi$ at the point $u$.
If $U$ is a Banach space, then we exploit the usual
notation $\langle\xi,u\rangle_{U^\ast,U}:=\xi(u)$.

\bigskip

\begin{example}\label{ex:banach}
Let $U$ and $V$ be Banach spaces and let $\tau_U$ and $\tau_V$ be the corresponding
weak topologies, i.e.
\begin{equation}
u_k\rightharpoonup u\quad\Leftrightarrow\quad
\langle\xi,u_k\rangle_{U^\ast,U}\to\langle\xi,u\rangle_{U^\ast,U}\;\forall\xi\in U^\ast.
\end{equation}
Then the similarity functional
\begin{equation} \label{eq:normex}
\cS(v_1,v_2):=\Vert v_1-v_2\Vert_V
\end{equation}
fulfills (iv) in Assumption~\ref{as:basic} with $s=1$. For a further discussion of this example
see Section~\ref{sec:source} below.
\end{example}

\bigskip

The next example, which is taken from \cite{Poe08}, shows that next to norms
also other similarity functionals are of interest.

\begin{example}
Let $(X,\rho)$ be a complete, separable metric space.
By $B(X)$ we denote the family of all Borel subsets of $X$, i.e. $B(X)$ is the
$\sigma$-algebra generated by the $\rho$-open sets in $X$, and by $P(X)$ we denote
the family of all Borel probability measures on $X$, i.e. $P(X)$ is the family of all measures
$\mu:B(X)\rightarrow[0,\infty)$ satisfying $\mu(X)=1$.
For $1\leq q<\infty$ we set
\begin{equation}
V:=\left\{\mu\in P(X):\int_X\rho(\mybullet,0)^q\diff\mu<\infty\right\}
\end{equation}
and as topology $\tau_V$ we choose the \emph{narrow topology} on $V$,
i.e. a series $(\mu_k)_{k\in\mathbb{N}}$ in $V$ converges to $\mu\in V$
with respect to $\tau_V$ if and only if
\begin{equation}
\int_X f\diff\mu_k\to\int_X f\diff\mu
\end{equation}
for all continuous and bounded real functions $f$ defined on $X$.
\par For defining the similarity functional $\cS$ we introduce the set
$\Gamma(\mu_1,\mu_2)\subseteq P(X\times X)$ (with $\mu_1,\mu_2\in V$) consisting of all
measures $\underline{\mu}\in P(X\times X)$ satisfying
$\underline{\mu}((\pi_i)^{-1}(A))=\mu_i(A)$ for all $A\in B(X)$ and $i=1,2$, where
$\pi_1(x_1,x_2):=x_1$ and $\pi_2(x_1,x_2):=x_2$.
Then the similarity functional $\cS$ defined by
\begin{equation}
\cS(\mu_1,\mu_2):=\left(\inf_{\underline{\mu}\in\Gamma(\mu_1,\mu_2)}
\int_{X\times X}\rho^q\diff\underline{\mu}\right)^q,\quad\mu_1,\mu_2\in V,
\end{equation}
is a metric, the \emph{Wasserstein metric}, on V and according to \cite{Poe08}
it satisfies (iv) in Assumption~\ref{as:basic}.
\par This similarity functional has been applied to flow,
mass transport, and image registration problems. For details on applications
and some references we refer to \cite{Poe08}.
\end{example}

\bigskip

In connection with the exponent $p$ in (\ref{eq:tikh}) from time to time we will make
use of the inequality
\begin{equation}
(a+b)^p\leq c_p(a^p+b^p)
\end{equation}
for $a\geq 0$ and $b\geq 0$ with
\begin{equation}\label{eq:cp}
c_p:=\begin{cases}1&\text{if }0<p<1,\\2^{p-1}&\text{if }p\geq 1.\end{cases}
\end{equation}

Under Assumption~\ref{as:basic} one can show that there exist
minimizers of the Tikhonov functional (\ref{eq:tikh}) for all $p>0$
and that these minimizers are stable with respect to perturbations
of the data $v^\delta$. The ideas of corresponding proofs given in Section~2 of \cite{HofKalPoeSch07}
and in \cite{Poe08} can be applied to our general setting, and we note that hence existence and stability
of minimizers can be ensured also in the case of exponents $0<p<1$ in Example~\ref{ex:banach}, for which assertions are up to now
missing in the literature.

To formulate assertions about convergence of a series of minimizers as $\delta$ tends
to zero we need the concept
of $\Omega$-minimizing solutions: An element $u^\dagger\in D$ is called
\emph{$\Omega$-minimizing solution} if
\begin{equation}
F(u^\dagger)=v^0\quad\text{and}\quad
\Omega(u^\dagger)=\inf\{\Omega(u):u\in D,\,F(u)=v^0\}.
\end{equation}
Under Assumption~\ref{as:basic} one can show that there exists an $\Omega$-minimizing solution.
The following theorem was proven in \cite{Poe08}.

\begin{theorem}
Assume that Assumption~\ref{as:basic} is satisfied. Let $(\delta_k)_{k\in\mathbb{N}}$ be
a sequence in $\mathbb{R}$ monotonically decreasing to zero, let
$\alpha:(0,\delta_1]\rightarrow(0,\infty)$ be a parameter choice with
$\alpha(\delta)\to 0$ and $\frac{\delta^p}{\alpha(\delta)}\to 0$ as $\delta\to 0$,
and let $\alpha_k:=\alpha(\delta_k)$ and $v_k:=v^{\delta_k}$.
Then every sequence $(u_k)_{k\in\mathbb{N}}$ in $U$ with
$u_k\in\argmin\{\cS(F(u),v_k)^p+\alpha\Omega(u):u\in D\}$ has a $\tau_U$-convergent
subsequence and the limit of each $\tau_U$-convergent subsequence is an $\Omega$-minimizing
solution. If the $\Omega$-minimizing solution is unique then $(u_k)$ converges
to this $\Omega$-minimizing solution.
\end{theorem}

To express convergence rates we use Bregman distances, which have
become quite popular in recent years for this purpose. In this
context, let
$$\tilde{u}\in D_B:=\{u\in U:\partial\Omega(u)\neq\emptyset\} \quad
\mbox{and} \quad \xi\in\partial\Omega(\tilde{u})\subseteq U^\ast,$$
where $\partial\Omega(u)$ denotes the subdifferential of $\Omega$ at $u$.
Then the functional
$$\cB_\xi(u,\tilde{u}):=\Omega(u)-\Omega(\tilde{u})-\xi(u-\tilde{u}),
\quad u\in U,$$ is called \emph{Bregman distance} with respect to
$\Omega$, $\tilde{u}$, and $\xi$. In the sequel we always assume
that there exists an $\Omega$-minimizing solution $u^\dagger\in
D_B$.

At the end of this section we want to mention the two basic concepts
occurring in the literature for proving convergence rates of
Tikhonov type regularization of ill-posed equations.
Classical source conditions, as, e.g., in a Banach space setting $\xi =
F^\prime(u^\dagger)^\ast\, \eta,\,$ $\eta\in V^\ast,$  and in Hilbert spaces
$u^\dagger = \varphi(F^\prime(u^\dagger)^\ast F^\prime(u^\dagger))\, w,$ $w\in U$,
i.e.~sourcewise representations of
an element $\xi$ of the subdifferential of $\Omega$ for an
$\Omega$-minimizing solution or of an $\Omega$-minimizing solution itself,
are the main ingredient for proving convergence rates.
%Classical source conditions, as, e.g., in a Banach space setting $\xi =
%F^\prime(u^\dagger)^\ast\, v,\,$ $v\in V^\ast,$  and in Hilbert spaces
%$\xi = \varphi(F^\prime(u^\dagger)^\ast F^\prime(u^\dagger))\, w,$ $w\in U$,
%i.e.~sourcewise representations of
%an element $\xi$ of the subdifferential of $\Omega$ for an
%$\Omega$-minimizing solution, are the main ingredient for proving
%convergence rates.
Such classical kinds of source conditions express
the smoothness of the solution with respect to the operator and they
alone are responsible for possible convergence rates in linear
ill-posed equation. If we are concerned with nonlinear operators $F$
then in addition we have to take into account structural conditions
which express the nonlinearity. For nonlinear ill-posed equations
classical source conditions and nonlinearity conditions together
control convergence rates. Their interplay, however, is rather
complicated.

Originally in~\cite{HofKalPoeSch07} (see also
\cite{SchGraGroHalLen09}) an extended concept of source conditions
was presented for obtaining convergence rates for the Banach space situation
of Example~\ref{ex:banach}. It
is based on variational inequalities, which have to hold on
appropriate level sets $M_\alpha(c)$ of the Tikhonov type functional
(\ref{eq:tikh}). In \cite{HeiHof09,Hof09a} an additional exponent
$\kappa \in (0,1]$ was motivated such that the variational
inequalities attain the form
\begin{equation} \label{eq:vinenorm}
\langle\xi,u^\dagger-u\rangle_{U^\ast,U} \le
\beta_1\cB_\xi(u,u^\dagger)+\beta_2
\|F(u)-F(u^\dagger)\|^\kappa_V.
\end{equation}
If such a variational inequality
holds, then a convergence rate
$\cB_\xi(u_{\alpha(\delta)}^\delta,u^\dagger)=\cO(\delta^\kappa)$ as
$\delta\to 0$ can immediately be derived without additional
knowledge when appropriate a priori parameter choices are used. Both
the classical source conditions and the structural conditions of
nonlinearity result into one parameter, namely the exponent $\kappa$
that alone controls the rate.

\section{Variational inequalities and convergence rates}\label{sec:vari}

In this section we extend the concept of variational inequalities
introduced in \cite{HofKalPoeSch07}. At first we state some simple
properties of the level sets defined in (\ref{eq:levelsets}).

\begin{proposition} \label{prop:levelsets}
Let $u^\dagger$ be an $\Omega$-minimizing solution and let
$\varrho>0$ be an arbitrary constant. Then for
$0<\alpha_1\leq\alpha_2$ we have
\begin{itemize}
\item[(i)]
$M_{\alpha_1}(\varrho)\supseteq M_{\alpha_2}(\varrho)$,
\item[(ii)]
$M_{\alpha_1}(\varrho\alpha_1)\subseteq M_{\alpha_2}(\varrho\alpha_2)$,
\item[(iii)]
$\bigcap\limits_{\alpha>0}M_\alpha(\varrho\alpha)=\{u\in D:F(u)=v^0,\,\Omega(u)\leq\varrho\}$.
\end{itemize}
\end{proposition}

\begin{proof}
Item (i) is trivial. Item (ii) follows from
\begin{align*}
\cT_{\alpha_2}^0(u)
&=\cS(F(u),v^0)^p+\alpha_1\Omega(u)-(\alpha_1-\alpha_2)\Omega(u)\\
&\leq\varrho\alpha_1-(\alpha_1-\alpha_2)\Omega(u)
=\varrho\alpha_2+(\alpha_1-\alpha_2)(\varrho-\Omega(u))
\leq\varrho\alpha_2
\end{align*}
for all $u\in M_{\alpha_1}(\varrho\alpha_1)$ and (iii) from
$\cS(F(u),v^0)^p\leq\alpha(\varrho-\Omega(u))$
for all $\alpha>0$ and for $u\in\bigcap_{\alpha>0}M_\alpha(\varrho\alpha)$.
\end{proof}

The next proposition shows the importance of the level sets
$M_\alpha(\varrho\alpha)$.

\begin{proposition}\label{th:param_choice}
Let $u^\dagger$ be an $\Omega$-minimizing solution, $\bar{\alpha}>0$, and
\begin{equation} \label{eq:rhodef}
\varrho\,>\,c_p\,s^p\,\Omega(u^\dagger)\,.
\end{equation}
 Further let $\delta\mapsto\alpha(\delta)$
be an a priori parameter choice satisfying
\begin{equation} \label{eq:apri}
\alpha(\delta)\to 0,\quad\frac{\delta^p}{\alpha(\delta)}\to 0\quad\text{as $\delta\to 0$}
\end{equation}
and let $u_{\alpha(\delta)}^\delta\in\argmin\{\cT_{\alpha(\delta)}^\delta(u):u\in D\}$
for $\delta>0$.
Then there exists a $\bar{\delta}>0$, such that
$u_{\alpha(\delta)}^\delta\in M_{\bar{\alpha}}(\varrho\bar{\alpha})$ holds for all
$\delta\in(0,\bar{\delta}]$.
\end{proposition}

\begin{proof}
Because $\alpha(\delta)\to 0$ and
$\frac{\delta^p}{\alpha(\delta)}\to 0$ as $\delta\to 0$ there exists
a $\bar{\delta}>0$ with $\alpha(\delta)\leq\bar{\alpha}$ and
$\frac{\delta^p}{\alpha(\delta)}\leq\frac{\varrho}{2c_ps^p}-\frac{1}{2}\Omega(u^\dagger)$
for all $\delta\in(0,\bar{\delta}]$. For the sake of brevity we
write $\alpha$ instead of $\alpha(\delta)$. For
$\delta\in(0,\bar{\delta}]$ we now have (in analogy to
\cite[p.\,5]{HeiHof09})
\begin{align*}
\cT_{\bar{\alpha}}^0(u_\alpha^\delta)
&\leq\bigl(s\cS(F(u_\alpha^\delta),v^\delta)+s\delta\bigr)^p+\bar{\alpha}\Omega(u_\alpha^\delta)\\
&\leq c_ps^p\bigl(\cS(F(u_\alpha^\delta),v^\delta)^p+\alpha\Omega(u_\alpha^\delta)
+(\bar{\alpha}-\alpha)\Omega(u_\alpha^\delta)+\delta^p\bigr)\\
&\leq c_ps^p\Bigl(2\delta^p+\alpha\Omega(u^\dagger)+\tfrac{\bar{\alpha}-\alpha}{\alpha}
\bigl(\cS(F(u_\alpha^\delta),v^\delta)^p+\alpha\Omega(u_\alpha^\delta)\bigr)\Bigr)\\
&=c_ps^p\bigl(\delta^p+\tfrac{\bar{\alpha}}{\alpha}\delta^p+\bar{\alpha}\Omega(u^\dagger)\bigr)
\leq c_ps^p\bar{\alpha}\bigl(2\tfrac{\delta^p}{\alpha}+\Omega(u^\dagger)\bigr)
\leq\varrho\bar{\alpha}.\qedhere
\end{align*}
\end{proof}

We now give the basic definition of a variational inequality in a
stronger sense.

\begin{definition}\label{def:vari}
An $\Omega$-minimizing solution $u^\dagger$ \emph{satisfies a
variational inequality}, if there exist a
$\xi\in\partial\Omega(u^\dagger)$ and constants
$\varrho$ fulfilling inequality (\ref{eq:rhodef}), $\bar{\alpha}>0$,
$\beta_1\in[0,1)$, $\beta_2\geq 0$, and $\kappa>0$, such that
\begin{equation}\label{eq:vari}
-\xi(u-u^\dagger)\leq\beta_1\cB_\xi(u,u^\dagger)+\beta_2\cS(F(u),F(u^\dagger))^\kappa
\end{equation}
holds for all $u\in M_{\bar{\alpha}}(\varrho\bar{\alpha})$.
\end{definition}

As one would expect, a variational inequality with $\kappa=\kappa_0$
implies a variational inequality with $\kappa=\kappa_1$ for each
$\kappa_1\in(0,\kappa_0)$. The only changing constant in Definition~\ref{def:vari}
is the factor $\beta_2=\beta_2(\kappa)$. This follows
immediately from
\begin{align*}
\beta_2(\kappa_0)\cS(F(u),F(u^\dagger))^{\kappa_0}
&=\beta_2(\kappa_0)\cS(F(u),F(u^\dagger))^{\kappa_0-\kappa_1}\cS(F(u),F(u^\dagger))^{\kappa_1}\\
&\leq\beta_2(\kappa_0)(\varrho\bar{\alpha})^{\frac{\kappa_0-\kappa_1}{p}}
\cS(F(u),F(u^\dagger))^{\kappa_1}
\end{align*}
because $u\in M_{\bar{\alpha}}(\varrho\bar{\alpha})$ implies
$\cS(F(u),F(u^\dagger))^p\leq\varrho\bar{\alpha}$.

For $\kappa=1$ and for the case of topological spaces with
general similarity functional $\cS$ Definition~\ref{def:vari} was introduced in \cite{Poe08}.
The definition was already presented earlier in \cite{HofKalPoeSch07} for the Banach space situation with norm as similarity
functional $\cS$. For that situation and $\kappa\in(0,1]$ this
variational inequality (\ref{eq:vari}) appeared also in \cite[proof of
Theorem~3.3]{HeiHof09}.

The connection between classical source conditions and variational inequalities
will be discussed in Section \ref{sec:source}.

We now give a first convergence rate result, which will be proven later in a more
general context.

\begin{theorem}\label{th:vari_rates}
Let $u^\dagger$ be an $\Omega$-minimizing solution which satisfies a
variational inequality in the sense of Definition~\ref{def:vari}
with $0<\kappa<p$ and let $\delta\mapsto\alpha(\delta)$ be an a priori
parameter choice with
$\underline{c}\delta^{p-\kappa}\leq\alpha(\delta)\leq\overline{c}\delta^{p-\kappa}$
for sufficiently small $\delta$ and constants $\underline{c}>0$,
$\overline{c}>0$. Then
\begin{equation}
\cB_\xi(u_{\alpha(\delta)}^\delta,u^\dagger)=\cO(\delta^\kappa)
\quad\text{as $\delta\to 0$}.
\end{equation}
\end{theorem}

\bigskip

Note that the a priori parameter choice in Theorem~\ref{th:vari_rates} restricts the
admissible values for the exponent $\kappa$ to the interval $(0,p)$. As we will see this restriction is due to the
 proof technique using Young's inequality. On the other hand, the following proposition provides an upper bound for $\kappa$ in a variational inequality
(\ref{eq:vari}). A special case of this proposition was also formulated in \cite[Proposition 4.3]{Hof09a}.

\begin{proposition}\label{th:kappa_bound}
Let $u^\dagger$ be an $\Omega$-minimizing solution which satisfies a
variational inequality in the sense of Definition~\ref{def:vari}. If
there exist a $q>0$, a $u\in U$ with $\xi(u)<0$, and a $t_0>0$,
such that $u^\dagger+tu\in M_{\bar{\alpha}}(\varrho\bar{\alpha})$
holds for all $t\in[0,t_0]$, the limits
$$L_\Omega:=\lim_{t\to+0}\frac{\Omega(u^\dagger+tu)-\Omega(u^\dagger)}{t},\quad
L_\cS:=\lim_{t\to+0}\frac{\cS(F(u^\dagger+tu),F(u^\dagger))^q}{t},$$
i.e. the directional derivatives in $u^\dagger$ in direction $u$ of $\Omega$
and $\cS(F(\mybullet),F(u^\dagger))^q$, exist, and $L_\Omega=\xi(u)$ holds,
then $\kappa\leq q$ must hold.
\end{proposition}

\begin{proof}
Let $\kappa>q$. For each $t\in(0,t_0]$ inequality (\ref{eq:vari}) then implies
$$-\xi(tu)\leq\beta_1\bigl(\Omega(u^\dagger+tu)-\Omega(u^\dagger)-\xi(tu)\bigr)
+\beta_2\cS(F(u^\dagger+tu),F(u^\dagger))^\kappa$$
and thus
$$-\xi(u)\leq\beta_1\underbrace{\left(\tfrac{\Omega(u^\dagger+tu)-\Omega(u^\dagger)}{t}-\xi(u)\right)}
_{\stackrel{t\to+0}{\longrightarrow}0}
+\beta_2\underbrace{\left(\tfrac{\cS(F(u^\dagger+tu),F(u^\dagger))^q}{t}\right)^{\frac{\kappa}{q}}}
_{\stackrel{t\to+0}{\longrightarrow}L_\cS^{\kappa/q}}\underbrace{t^{\frac{\kappa}{q}-1}}
_{\stackrel{t\to+0}{\longrightarrow}0}.$$
Passage to the limit $t\to+0$ gives $\xi(u)\geq 0$, which is a contradiction to $\xi(u)<0$.
\end{proof}

\begin{remark} \label{rem:rem1}
Under the standing assumptions of this paper on $F,\,D(F),\, \Omega,$ and $u^\dagger$
one can easily show that for Banach spaces $U$ and $V$ and $\cS(v_1,v_2):=\Vert v_1-v_2\Vert_V$ the Proposition~\ref{th:kappa_bound} applies
for $q=1$ when $F$ and $\Omega$ are G\^ateaux differentiable in $u^\dagger$.
So in this case only variational inequalities (\ref{eq:vinenorm}) with $\kappa\leq 1$ can be satisfied if the singular case $\xi=\Omega^\prime(u^\dagger)=0$ is excluded.
\end{remark}

\section{Approximate variational inequalities}\label{sec:appvari}

The aim of this section is to formulate convergence rates results
without assuming that a variational inequality is satisfied. As in
the method of approximate source conditions (see \cite{DHY07} and
\cite{HeiHof09}) we use distance functions
$d:[0,\infty)\rightarrow[0,\infty)$ measuring the violation of a
prescribed benchmark condition. However, here we have a variational
inequality (\ref{eq:vari}) as benchmark condition and the distance
functions are defined in a completely different manner.

If a benchmark inequality of type (\ref{eq:vari}) is not satisfied
then there exists at least one $u\in
M_{\bar{\alpha}}(\varrho\bar{\alpha})$ with
$$-\xi(u-u^\dagger)
>\beta_1\cB_\xi(u,u^\dagger)+\beta_2\cS(F(u),F(u^\dagger))^\kappa.$$
Thus the ``maximum violation'' of a variational inequality (\ref{eq:vari})
may be expressed by
\begin{equation}
\sup_{u\in M_{\bar{\alpha}}(\varrho\bar{\alpha})}\bigl(-\xi(u-u^\dagger)
-\beta_1\cB_\xi(u,u^\dagger)-\beta_2\cS(F(u),F(u^\dagger))^\kappa\bigr).
\end{equation}
The question whether the satisfaction of the benchmark inequality can
be forced by increasing the factor $\beta_2$ leads to the definition
of an approximate variational inequality.

\begin{definition}\label{def:appvari}
An $\Omega$-minimizing solution \emph{satisfies an approximate
variational inequality} (approximate inequality for short) if there
exist a $\xi\in\partial\Omega(u^\dagger)$ and constants
$\varrho$ fulfilling (\ref{eq:rhodef}), $\bar{\alpha}>0$,
$\beta_1\in[0,1)$, $\beta_2\geq 0$, $\gamma\geq 0$ and $\kappa>0$,
such that the function $d:[0,\infty)\rightarrow\mathbb{R}$ defined
by
$$d(r):=-\min_{u\in M_{\bar{\alpha}}(\varrho\bar{\alpha})}\bigl(\xi(u-u^\dagger)
+\beta_1\cB_\xi(u,u^\dagger)+\beta_2r^\gamma\cS(F(u),F(u^\dagger))^\kappa\bigr)$$
satisfies $d(r)\to 0$ as $r\to\infty$.
\end{definition}

The constant $\gamma$ in Definition~\ref{def:appvari} seems to be not necessary,
but it will turn out that it explicitly occurs in the formulation of convergence rates.

At first we prove some basic properties of the distance function $d$.

\begin{proposition}\label{th:dprop}
Let $u^\dagger$ be an $\Omega$-minimizing solution which satisfies
an approximate inequality in the sense of Definition~\ref{def:appvari}. Then we have:
\begin{itemize}
\item[(i)]
$0\leq d(r)<\infty$ holds for all $r\geq 0$.
\item[(ii)]
The minimum in the definition of $d$ is attained.
\item[(iii)]
$d$ is continuous.
\item[(iv)]
$d$ is monotonically decreasing.
\item[(v)]
If $d(r)>0$ holds for all $r\geq 0$, then $d$ is strictly monotonically decreasing.
\end{itemize}
\end{proposition}

\begin{proof}
\mbox{}
\begin{itemize}
\item[(i)]
Because $\cT_{\bar{\alpha}}^0(u^\dagger)=\bar{\alpha}\Omega(u^\dagger)\leq\varrho\bar{\alpha}$
we have $u^\dagger\in M_{\bar{\alpha}}(\varrho\bar{\alpha})$ and therefore $d(r)\geq 0$.
For $r\geq 0$ from $\Omega(u)\leq\varrho$ and
$\cS(F(u),F(u^\dagger))\leq(\varrho\bar{\alpha})^{\frac{1}{p}}$
for $u\in M_{\bar{\alpha}}(\varrho\bar{\alpha})$ we get the estimate
\begin{equation}
d(r)\leq C(r,\Omega(u^\dagger))
+\sup_{u\in M_{\bar{\alpha}}(\varrho\bar{\alpha})}\vert\xi(u)\vert
\end{equation}
with a constant $C<\infty$ depending on $r$ and $\Omega(u^\dagger)$ only.
Assume there exists a sequence $(u_k)_{k\in\mathbb{N}}$ in $M_{\bar{\alpha}}(\varrho\bar{\alpha})$
with $\vert\xi(u_k)\vert\to\infty$. Then from Assumption~\ref{as:basic} (vii)
the existence of a $\tau_U$-convergent subsequence $(u_{k_l})_{l\in\mathbb{N}}$ follows;
let $\tilde{u}\in U$ be its limit.
The continuity of $\xi$ implies
$\vert\xi(u_{k_l})\vert\to\vert\xi(\tilde{u})\vert$
and therefore the boundedness of the sequence $(\vert\xi(u_{k_l})\vert)$.
This contradicts $\vert\xi(u_k)\vert\to\infty$. Thus
$$\sup_{u\in M_{\bar{\alpha}}(\varrho\bar{\alpha})}\vert\xi(u)\vert<\infty,$$
i.e. $d(r)<\infty$.
\item[(ii)]
We define $g_r:U\rightarrow\mathbb{R}\cup\{+\infty\}$ by
\begin{equation}
g_r(u):=\xi(u-u^\dagger)+\beta_1\cB_\xi(u,u^\dagger)
+\beta_2r^\gamma\cS(F(u),F(u^\dagger))^\kappa.
\end{equation}
The continuity of $\xi$ and $F$ and the lower semi-continuity of $\Omega$ and $\cS$
together imply the lower semi-continuity of $g_r$.
Now let $(u_k)_{k\in\mathbb{N}}$ be a sequence in $M_{\bar{\alpha}}(\varrho\bar{\alpha})$
satisfying $g_r(u_k)\to\inf_{u\in M_{\bar{\alpha}}(\varrho\bar{\alpha})}g_r(u)$.
Then there exists a $\tau_U$-convergent subsequence $(u_{k_l})_{l\in\mathbb{N}}$ with
limit $\tilde{u}\in U$, especially it holds $\tilde{u}\in M_{\bar{\alpha}}(\varrho\bar{\alpha})$
(because $\cT_{\bar{\alpha}}^0$ is lower semi-continuous), and
$$g_r(\tilde{u})\leq\liminf_{l\to\infty}g_r(u_{k_l})=\lim_{l\to\infty}g_r(u_{k_l})
=\inf_{u\in M_{\bar{\alpha}}(\varrho\bar{\alpha})}g_r(u)$$
holds. Thus, $g_r(\tilde{u})=\inf_{u\in M_{\bar{\alpha}}(\varrho\bar{\alpha})}g_r(u)$.
\item[(iii)]
For $r\geq 0$ let $g_r$ be defined as in the proof of item (ii) and let
$u_r\in M_{\bar{\alpha}}(\varrho\bar{\alpha})$ be a minimizer of $g_r$.
Then for all $r,s\geq 0$ we have
$$d(r)-d(s)=\underbrace{\min g_s}_{\leq g_s(u_r)}-\underbrace{\min g_r}_{=g_r(u_r)}
\leq\beta_2(s^\gamma-r^\gamma)
{\underbrace{\cS(F(u_r),F(u^\dagger))}_{\leq(\varrho\bar{\alpha})^{\frac{1}{p}}}}^\kappa$$
and
$$-(d(r)-d(s))=\underbrace{\min g_r}_{\leq g_r(u_s)}-\underbrace{\min g_s}_{=g_s(u_s)}
\leq\beta_2(r^\gamma-s^\gamma)
{\underbrace{\cS(F(u_s),F(u^\dagger))}_{\leq(\varrho\bar{\alpha})^{\frac{1}{p}}}}^\kappa,$$
i.e.
\begin{equation}
\vert d(r)-d(s)\vert\leq\beta_2(\varrho\bar{\alpha})^{\frac{\kappa}{p}}
\vert r^\gamma-s^\gamma\vert,
\end{equation}
implying the continuity of $d$.
\item[(iv)]
The assertion follows directly from the definition of $d$.
\item[(v)]
If $\beta_2=0$ or $\gamma=0$ would hold, then $d$ would be constant. But this is
not possible because $d(r)>0$ for all $r\geq 0$ and $d(r)\to 0$ as $r\to\infty$.
Thus $\beta_2>0$ and $\gamma>0$ hold.
\par We assume that there is an $r\geq 0$ for which $g_r$ (set as in the proof of (ii))
has a minimizer $\tilde{u}\in M_{\bar{\alpha}}(\varrho\bar{\alpha})$
satisfying $F(\tilde{u})=v^0$. Then for each $s\geq 0$ we get
\begin{align*}
\lefteqn{\xi(\tilde{u}-u^\dagger)+\beta_1\cB_\xi(\tilde{u},u^\dagger)}\\
&\quad=\xi(\tilde{u}-u^\dagger)+\beta_1\cB_\xi(\tilde{u},u^\dagger)
+\beta_2s^\gamma\cS(F(\tilde{u}),F(u^\dagger))^\kappa\\
&\quad\geq\min_{u\in M_{\bar{\alpha}}(\varrho\bar{\alpha})}
\bigl(\xi(u-u^\dagger)+\beta_1\cB_\xi(u,u^\dagger)
+\beta_2s^\gamma\cS(F(u),F(u^\dagger))^\kappa\bigr)\\
&\quad=-d(s)
\end{align*}
and thus $d(s)\to 0$ as $s\to\infty$ implies
$\xi(\tilde{u}-u^\dagger)+\beta_1\cB_\xi(\tilde{u},u^\dagger)\geq 0$.
But this contradicts
\begin{align*}
\lefteqn{\xi(\tilde{u}-u^\dagger)+\beta_1\cB_\xi(\tilde{u},u^\dagger)}\\
&\quad=\xi(\tilde{u}-u^\dagger)+\beta_1\cB_\xi(\tilde{u},u^\dagger)
+\beta_2r^\gamma\cS(F(\tilde{u}),F(u^\dagger))^\kappa\\
&\quad=\min_{u\in M_{\bar{\alpha}}(\varrho\bar{\alpha})}
\bigl(\xi(u-u^\dagger)+\beta_1\cB_\xi(u,u^\dagger)
+\beta_2r^\gamma\cS(F(u),F(u^\dagger))^\kappa\bigr)\\
&\quad=-d(r)<0.
\end{align*}
So for each $r\geq 0$ each minimizer $\tilde{u}\in M_{\bar{\alpha}}(\varrho\bar{\alpha})$
of $g_r$ satisfies the inequality $\cS(F(\tilde{u}),F(u^\dagger))>0$.
Now for $0\leq s<r$ we have
\begin{align*}
d(r)&=-\min_{u\in M_{\bar{\alpha}}(\varrho\bar{\alpha})}g_r(u)
=-g_r(\tilde{u})\\
&=-g_s(\tilde{u})-\underbrace{\beta_2(r^\gamma-s^\gamma)}_{>0}
\underbrace{\cS(F(\tilde{u}),F(u^\dagger))^\kappa}_{>0}\\
&<-g_s(\tilde{u})\leq-\min_{u\in M_{\bar{\alpha}}(\varrho\bar{\alpha})}g_s(u)
=d(s),
\end{align*}
i.e. $d$ ist strictly monotonically decreasing.\qedhere
\end{itemize}
\end{proof}

Obviously an $\Omega$-minimizing solution satisfies a variational
inequality in the sense of Definition~\ref{def:vari} if and only if
it satisfies an approximate inequality in the sense of Definition~\ref{def:appvari}
and there exists an $r_0\geq 0$ with $d(r_0)=0$.

If $u^\dagger$ satifies an approximate inequality with constant
$\bar{\alpha}=\alpha_0$ then it satisfies an approximate inequality
with $\bar{\alpha}=\alpha_1$ for all $\alpha_1\in(0,\alpha_0]$ and
with the same other constants. Later we will see that the constant
$\bar{\alpha}$ from Definition~\ref{def:appvari} does not appear
explicitly in the formulation of convergence rates. So for the sake
of plausibility of Definition~\ref{def:appvari} the distance
function $d$ should be independent of $\bar{\alpha}$. The next two
propositions give some insight into this problem.

\begin{proposition}\label{th:concent1}
Let $u^\dagger$ be an $\Omega$-minimizing solution which satisfies
an approximate inequality in the sense of Definition~\ref{def:appvari}.
Further let $(r_k)_{k\in\mathbb{N}}$ be a
sequence in $(0,\infty)$ with $r_k\to\infty$ and let
$(u_k)_{k\in\mathbb{N}}$ be a sequence of elements $u_k\in
M_{\bar{\alpha}}(\varrho\bar{\alpha})$ which realize the minimum in
the definition of $d$, such that $u_k\rightharpoonup\tilde{u}$ holds
for some $\tilde{u}\in D$. Then it follows
$$F(\tilde{u})=v^0,\quad\Omega(\tilde{u})\leq\varrho,\quad\text{and}\quad
\xi(\tilde{u}-u^\dagger)=\tfrac{-\beta_1}{1-\beta_1}(\Omega(\tilde{u})
-\Omega(u^\dagger)).$$
\end{proposition}

\begin{proof}
The definitions of $u_k$ and $d(r_k)$ imply
$$-\beta_2r_k^\gamma\cS(F(u_k),F(u^\dagger))^\kappa=d(r_k)
+\xi(u_k-u^\dagger)
+\beta_1\cB_\xi(u_k,u^\dagger).$$
>From the continuity of $\xi$ and the lower semi-continuity of $\Omega$ for $\varepsilon>0$
and suffiently large $k\in\mathbb{N}$ it follows
$$-\beta_2r_k^\gamma\cS(F(u_k),F(u^\dagger))^\kappa\geq d(r_k)
+\xi(\tilde{u}-u^\dagger)
+\beta_1\cB_\xi(\tilde{u},u^\dagger)-\varepsilon$$
and therefore
$$\cS(F(u_k),F(u^\dagger))^\kappa\leq\underbrace{\tfrac{-1}{\beta_2r_k^\gamma}}_{\to 0}
\bigl(\underbrace{d(r_k)}_{\to 0}+\xi(\tilde{u}-u^\dagger)
+\beta_1\cB_\xi(\tilde{u},u^\dagger)-\varepsilon\bigr).$$
Passage to the limit $k\to\infty$ gives
$\cS(F(u_k),F(u^\dagger))^\kappa\to 0$
and with Assumption~\ref{as:basic} (iv)(d) this implies $F(u_k)\rightharpoonup v^0$.
On the other hand Assumption~\ref{as:basic} (i) implies $F(u_k)\rightharpoonup F(\tilde{u})$
and therefore $F(\tilde{u})=v^0$ holds.
\par The second assertion follows from
$$\Omega(\tilde{u})\leq\liminf_{k\to\infty}\Omega(u_k)\leq\liminf_{k\to\infty}\tfrac{1}{\bar{\alpha}}
\cT_{\bar{\alpha}}^0(u_k)\leq\varrho.$$
\par To prove the third and last assertion we first observe
\begin{align*}
\lefteqn{-\xi(\tilde{u}-u^\dagger)-\beta_1\cB_\xi(\tilde{u},u^\dagger)}\\
&\quad=-\xi(\tilde{u}-u^\dagger)-\beta_1\cB_\xi(\tilde{u},u^\dagger)
-\beta_2r_k^\gamma\cS(F(\tilde{u}),F(u^\dagger))^\kappa\\
&\quad\leq d(r_k)\to 0,
\end{align*}
which gives
\begin{equation}\label{eq:misc01}
-\xi(\tilde{u}-u^\dagger)-\beta_1\cB_\xi(\tilde{u},u^\dagger)\leq 0.
\end{equation}
For $\varepsilon>0$ and $k\in\mathbb{N}$ sufficiently large the continuity of $\xi$ and
the lower semicontinuity of $\Omega$ imply
\begin{align*}
0&\geq-\beta_2r_k^\gamma\cS(F(u_k),F(u^\dagger))^\kappa
=d(r_k)+\xi(u_k-u^\dagger)+\beta_1\cB_\xi(u_k,u^\dagger)\\
&\geq\underbrace{d(r_k)}_{\to 0}+\xi(\tilde{u}-u^\dagger)
+\beta_1\cB_\xi(\tilde{u},u^\dagger)-\varepsilon.
\end{align*}
By passage to the limit $k\to\infty$ we get
$\xi(\tilde{u}-u^\dagger)
+\beta_1\cB_\xi(\tilde{u},u^\dagger)\leq\varepsilon$
and from the arbitrarity of $\varepsilon$ it follows
\begin{equation}\label{eq:misc02}
-\xi(\tilde{u}-u^\dagger)-\beta_1\cB_\xi(\tilde{u},u^\dagger)\geq 0.
\end{equation}
Inequalities (\ref{eq:misc01}) and (\ref{eq:misc02}) together imply
$$-\xi(\tilde{u}-u^\dagger)=\beta_1\cB_\xi(\tilde{u},u^\dagger)$$
and substituting the Bregman distance by its definition gives the assertion.
\end{proof}

\begin{proposition}\label{th:concent2}
Let $u^\dagger$ be an $\Omega$-minimizing solution which satisfies
an approximate inequality in the sense of Definition~\ref{def:appvari}
and let $d_\alpha$ for $\alpha\in(0,\bar{\alpha}]$
be the function defined in analogy to $d$ with $\bar{\alpha}$
replaced by $\alpha$. If there exists no $u\in U$ with $F(u)=v^0$,
$\cR(u)=\varrho$ and
$\xi(u-u^\dagger)=\frac{-\beta_1}{1-\beta_1}(\varrho-\Omega(u^\dagger))$
then the following assertions are true:
\begin{itemize}
\item[(i)]
For all $\alpha\in(0,\bar{\alpha}]$ there exists an $r_\alpha\geq 0$, such that
$d(r)=d_\alpha(r)$ holds for all $r\geq r_\alpha$.
\item[(ii)]
For all $\alpha\in(0,\bar{\alpha}]$ there exists an $r_\alpha\geq 0$, such that for all
$r\geq r_\alpha$ all elements of $M_{\bar{\alpha}}(\varrho\bar{\alpha})$ which realize
the minimum in the definition of $d(r)$ lie in $M_\alpha(\varrho\alpha)$.
\end{itemize}
\end{proposition}

\begin{proof}
Assertion (i) is a direct consequence of (ii). We give an indirect proof of assertion (ii).
We assume that there exist an $\alpha\in(0,\bar{\alpha}]$ and a sequence
$(r_k)_{k\in\mathbb{N}}$ in $(0,\infty)$ with $r_k\to\infty$, such that for each $r_k$
there exists  an element $u_k\in M_{\bar{\alpha}}(\varrho\bar{\alpha})$ which
realizes the minimum in the definition of $d(r_k)$ and which satisfies
$u_k\notin M_\alpha(\varrho\alpha)$.
Because of Assumption~\ref{as:basic} (vii) and the lower semi-continuity of $\cT_{\bar{\alpha}}^0$
the sequence $(u_k)_{k\in\mathbb{N}}$ has a convergent subsequence, which we again denote
by $(u_k)_{k\in\mathbb{N}}$, with limit $\tilde{u}\in M_{\bar{\alpha}}(\varrho\bar{\alpha})$.
\par Proposition \ref{th:concent1} now implies
\begin{equation}\label{eq:misc03}
F(\tilde{u})=v^0,\quad\Omega(\tilde{u})\leq\varrho\quad\text{and}\quad
\xi(\tilde{u}-u^\dagger)
=\tfrac{-\beta_1}{1-\beta_1}(\Omega(\tilde{u})-\Omega(u^\dagger)).
\end{equation}
>From $u_k\notin M_\alpha(\varrho\alpha)$ in addition it follows
$$\Omega(u_k)>\varrho-\tfrac{1}{\alpha}\cS(F(u_k),v^0)^p$$
for all $k\in\mathbb{N}$ and thus $\cS(F(u_k),F(u^\dagger))\to 0$
(c.f. proof of Proposition~\ref{th:concent1}) implies $\Omega(u_k)>\varrho-\varepsilon$
for $\varepsilon>0$ and sufficiently large $k\in\mathbb{N}$.
Together with $\Omega(u_k)\leq\varrho$ this gives $\Omega(u_k)\to\varrho$.
Therefore from
\begin{align*}
0&\geq-\beta_2r_k^\gamma\cS(F(u_k),F(u^\dagger))^\kappa\\
&=\underbrace{d(r_k)}_{\to 0}
+\xi(u_k-u^\dagger)
+\beta_1\bigl(\Omega(u_k)-\Omega(u^\dagger)-\xi(u_k-u^\dagger)\bigr)
\end{align*}
by passage to the limit we conclude
$$0\geq(1-\beta_1)\xi(\tilde{u}-u^\dagger)
+\beta_1(\varrho-\Omega(u^\dagger))$$
and together with (\ref{eq:misc03}) we get
\begin{align*}
\tfrac{-\beta_1}{1-\beta_1}(\Omega(\tilde{u})-\Omega(u^\dagger))
&=\xi(\tilde{u}-u^\dagger)
\leq\tfrac{-\beta_1}{1-\beta_1}(\varrho-\Omega(u^\dagger))\\
&\leq\tfrac{-\beta_1}{1-\beta_1}(\Omega(\tilde{u})-\Omega(u^\dagger)),
\end{align*}
i.e. especially it holds $\cR(\tilde{u})=\varrho$.
Substituting this equality into (\ref{eq:misc03}) gives a contradiction to the
assumptions of the proposition.
\end{proof}

The following Lemma prepares the main theorem of this paper.

\begin{lemma}\label{th:rates_lemma}
Let $u^\dagger$ be an $\Omega$-minimizing solution which satisfies
an approximate inequality in the sense of Definition~\ref{def:appvari} with $0<\kappa<p$. Further let
$\alpha\mapsto\alpha(\delta)$ be a parameter choice fulfilling the condition (\ref{eq:apri}) from
Proposition~\ref{th:param_choice} and let $\bar{\delta}$ be the
corresponding constant from that proposition. Then
there exist constants $K_1>0$, $K_2>0$, and $K_3>0$, such that
\begin{equation}
\cB_\xi(u_{\alpha(\delta)}^\delta,u^\dagger)\leq K_1\tfrac{\delta^p}{\alpha(\delta)}
+K_2\alpha(\delta)^{\frac{\kappa}{p-\kappa}}r^{\frac{\gamma p}{p-\kappa}}+K_3 d(r)
\end{equation}
holds for all $r\geq 0$ and all $\delta\in(0,\bar{\delta}]$.
\end{lemma}

\begin{proof}
For the sake of brevity we write $\alpha$ instead of $\alpha(\delta)$.
Proposition \ref{th:param_choice} and the definition of $d(r)$ give us the
inequality
\begin{equation}
-\xi(u_\alpha^\delta-u^\dagger)
\leq\beta_1\cB_\xi(u_\alpha^\delta,u^\dagger)
+\beta_2r^\gamma\cS(F(u_\alpha^\delta),F(u^\dagger))^\kappa+d(r)
\end{equation}
for sufficiently small $\delta$.
>From this we get
\begin{align*}
\alpha\cB_\xi(u_\alpha^\delta,u^\dagger)
&=\cS(F(u_\alpha^\delta),v^\delta)^p+\alpha\Omega(u_\alpha^\delta)-\alpha\Omega(u^\dagger)
-\alpha\xi(u_\alpha^\delta-u^\dagger)-\cS(F(u_\alpha^\delta),v^\delta)^p\\
&\leq\delta^p-\alpha\xi(u_\alpha^\delta-u^\dagger)
-\cS(F(u_\alpha^\delta),v^\delta)^p\\
&\leq\delta^p+\alpha\beta_1\cB_\xi(u_\alpha^\delta,u^\dagger)
+\alpha\beta_2r^\gamma\cS(F(u_\alpha^\delta),F(u^\dagger))^\kappa+\alpha d(r)
-\cS(F(u_\alpha^\delta),v^\delta)^p\\
&\leq\delta^p+\alpha\beta_1\cB_\xi(u_\alpha^\delta,u^\dagger)
+\alpha\beta_2r^\gamma s^\kappa c_\kappa\bigl(\cS(F(u_\alpha^\delta),v^\delta)^\kappa
+\delta^\kappa\bigr)\\
&\quad\,+\alpha d(r)-\cS(F(u_\alpha^\delta),v^\delta)^p,
\end{align*}
where $c_\kappa$ in analogy to $c_p$ is given by
$$c_\kappa:=\begin{cases}1&\text{if }0<\kappa<1,\\
2^{\kappa-1}&\text{if }\kappa\geq 1.\end{cases}$$
Thus we have
\begin{align} \label{eq:pgleichkappa}
\cB_\xi(u_\alpha^\delta,u^\dagger)
&\leq\tfrac{1}{\alpha(1-\beta_1)}\Bigl(2\delta^p
+\alpha c_\kappa\beta_2r^\gamma s^\kappa\delta^\kappa-\delta^p
+\alpha c_\kappa\beta_2r^\gamma s^\kappa\cS(F(u_\alpha^\delta),v^\delta)^\kappa\Bigr. \nonumber \\
&\quad\qquad\qquad\;\Bigl.-\cS(F(u_\alpha^\delta),v^\delta)^p+\alpha d(r)\Bigr).
\end{align}

Now we apply the inequality
\begin{equation}
ab-\varepsilon a^{p_1}\leq\frac{b^{p_2}}{(\varepsilon p_1)^{p_2/p_1}p_2},
\end{equation}
where $a,b\geq 0$, $\varepsilon>0$, $p_1,p_2>1$ and $\frac{1}{p_1}+\frac{1}{p_2}=1$
have to hold, once with
$$a:=\delta^\kappa,\quad b:=\alpha c_\kappa\beta_2 r^\gamma s^\kappa,\quad\varepsilon:=1,
\quad p_1:=\tfrac{p}{\kappa},\quad p_2:=\tfrac{p}{p-\kappa}$$
and once with $\cS(F(u_\alpha^\delta),v^\delta)$ instead of $\delta$.
We get
\begin{align*}
\cB_\xi(u_\alpha^\delta,u^\dagger)
&\leq\tfrac{1}{\alpha(1-\beta_1)}\left(2\delta^p
+2(c_\kappa\beta_2s^\kappa)^{\tfrac{p}{p-\kappa}}
\bigl(\tfrac{\kappa}{p}\bigr)^{\tfrac{\kappa}{p-\kappa}}\tfrac{p-\kappa}{p}
\alpha^{\tfrac{p}{p-\kappa}}r^{\tfrac{\gamma p}{p-\kappa}}
+\alpha d(r)\right)\\
&=\underbrace{\tfrac{2}{1-\beta_1}}_{=:K_1}\tfrac{\delta^p}{\alpha}
+\underbrace{2(c_\kappa\beta_2s^\kappa)^{\frac{p}{p-\kappa}}
\bigl(\tfrac{\kappa}{p}\bigr)^{\tfrac{\kappa}{p-\kappa}}\tfrac{p-\kappa}{p(1-\beta_1)}}_{=:K_2}
\alpha^{\tfrac{\kappa}{p-\kappa}}r^{\tfrac{\gamma p}{p-\kappa}}
+\underbrace{\tfrac{1}{1-\beta_1}}_{=:K_3}d(r).
\end{align*}
\end{proof}

Now we can prove the convergence rate theorem from Section \ref{sec:vari}.

\begin{proof}[Proof of Theorem \ref{th:vari_rates}]
Because $u^\dagger$ satisfies a variational inequality it also
satisfies an approximate inequality with a distance function $d$ for
which there exists an $r_0\geq 0$ with $d(r)=0$ for all $r\geq r_0$.
So the assertion follows immediately from Lemma~\ref{th:rates_lemma}
with $r:=r_0$.
\end{proof}

\begin{theorem}\label{th:appvari_rates}
Let $u^\dagger$ be an $\Omega$-minimizing solution which satisfies for some $0<\kappa<p$
an approximate inequality in the sense of Definition~\ref{def:appvari} with $d(r)>0$ for all $r\geq 0$.
For $r>0$ we define
\begin{equation}
\Psi(r):=d(r)^{\frac{p-\kappa}{\kappa}}r^{-\frac{\gamma p}{\kappa}}\quad\text{and}\quad
\Phi(r):=d(r)^{\frac{1}{\kappa}}r^{-\frac{\gamma}{\kappa}}.
\end{equation}
Further let $\alpha\mapsto\alpha(\delta)$ be a parameter choice with
$\delta^p=\alpha(\delta)d\bigl(\Psi^{-1}(\alpha(\delta))\bigr)$ for sufficiently
small $\delta>0$. Then
\begin{equation}
\cB_\xi(u_{\alpha(\delta)}^\delta,u^\dagger)=\cO\bigl(d(\Phi^{-1}(\delta))\bigr)
\quad\text{as $\delta\to 0$}
\end{equation}
holds.
\end{theorem}

\begin{proof}
For the sake of brevity we write $\alpha(\delta)$ instead of $\alpha$.
Because $d$ is strictly monotonically decreasing $\Psi$ and $\Phi$ are strictly monotonically
decreasing, too. Thus the inverse functions $\Psi^{-1}$ and $\Phi^{-1}$ exist
and are strictly monotonically decreasing.
\par Lemma \ref{th:rates_lemma} with $r:=\Psi^{-1}(\alpha)$, i.e.
$$\alpha^{\frac{\kappa}{p-\kappa}}r^{\frac{\gamma p}{p-\gamma}}
=\Psi(r)^{\frac{\kappa}{p-\kappa}}r^{\frac{\gamma p}{p-\gamma}}=d(r),$$
implies
$$\cB_\xi(u_\alpha^\delta,u^\dagger)\leq K_1\tfrac{\delta^p}{\alpha}
+(K_2+K_3)d(\Psi^{-1}(\alpha))
=(K_1+K_2+K_3)d(\Psi^{-1}(\alpha))$$
for sufficiently small $\delta\leq\bar{\delta}$ and from
\begin{align*}
\Phi(\Psi^{-1}(\alpha))
&=d(\Psi^{-1}(\alpha))^{\frac{1}{\kappa}}\Psi^{-1}(\alpha)^{-\frac{\gamma}{\kappa}}
=\bigl(\tfrac{\delta^p}{\alpha}\bigr)^{\frac{1}{\kappa}}
\Psi^{-1}(\alpha)^{-\frac{\gamma}{\kappa}}\\
&=\delta^{\frac{p}{\kappa}}{\underbrace{\left(\alpha^{\frac{\kappa}{p-\kappa}}
\Psi^{-1}(\alpha)^{\frac{\gamma p}{p-\kappa}}\right)}_{=d(\Psi^{-1}(\alpha))
=\frac{\delta^p}{\alpha}}}^{\frac{\kappa-p}{\kappa p}}
\alpha^{\frac{1}{p}-\frac{1}{\kappa}}
=\delta
\end{align*}
we conclude $\Psi^{-1}(\alpha)=\Phi^{-1}(\delta)$, which proves the assertion.
\end{proof}

\begin{remark}
If instead of $d$ only a strictly monotonically decreasing majorant $\bar{d}$ of $d$
is available then Lemma~\ref{th:rates_lemma} and Theorem~\ref{th:appvari_rates} also hold with
$d$ replaced by $\bar{d}$.
\end{remark}

The following propositions give some further insight into the convergence rates
results of this paper.

\begin{proposition}
With the notation of Theorem~\ref{th:appvari_rates} it holds
\begin{equation}
d(\Phi^{-1}(\delta))=\co(\delta^\kappa)\quad\text{as $\delta\to 0$}.
\end{equation}
\end{proposition}

\begin{proof}
With $r:=\Phi^{-1}(\delta)$, i.e. $\delta=\Phi(r)$, we have
\begin{equation}
\frac{\delta^\kappa}{d(\Phi^{-1}(\delta))}=\frac{\Phi(r)^\kappa}{d(r)}
=\frac{d(r)r^{-\gamma}}{d(r)}=r^{-\gamma}=\Phi^{-1}(\delta)^{-\gamma}.
\end{equation}
>From $\gamma>0$ and $\Phi^{-1}(\delta)\to\infty$ as $\delta\to 0$ we conclude
$\Phi^{-1}(\delta)^{-\gamma}\to 0$ as $\delta\to 0$ and therefore the assertion
follows.
\end{proof}

\begin{proposition}\label{th:d_at_infty}
With the notation of Theorem~\ref{th:appvari_rates}
\begin{equation}
d(r)=\cO(\tilde{d}(r))\quad\text{as $r\to\infty$}
\end{equation}
implies
\begin{equation}
d(\Phi^{-1}(\delta))=\cO\bigl(\tilde{d}(\tilde{\Phi}^{-1}(\delta))\bigr)
\quad\text{as $\delta\to 0$},
\end{equation}
where $\tilde{\Phi}$ is defined in the same way as $\Phi$ with $d$ replaced by
$\tilde{d}$.
\end{proposition}

\begin{proof}
>From $d(r)=\cO(\tilde{d}(r))$ as $r\to\infty$ it follows $\Phi(r)=\cO(\tilde{\Phi}(r))$
as $r\to\infty$ and this implies $d(\Phi^{-1}(\delta))=\cO(\tilde{\Phi}^{-1}(\delta))$
as $\delta\to 0$.
With $r:=\Phi^{-1}(\delta)$, i.e. $\delta=\Phi(r)$, for sufficiently small $\delta>0$
and a constant $c>0$ we get
\begin{align*}
d(\Phi^{-1}(\delta))
&=\delta^\kappa\frac{d(\Phi^{-1}(\delta))}{\delta^\kappa}
=\Phi(r)^\kappa\frac{d(r)}{d(r)r^{-\gamma}}\\
&=\delta^\kappa\Phi^{-1}(\delta)^\gamma
\leq c\delta^\kappa\tilde{\Phi}^{-1}(\delta)^\gamma=c\tilde{d}(\tilde{\Phi}^{-1}(\delta)).
\end{align*}
\end{proof}

Proposition \ref{th:kappa_bound} told us that under weak assumptions there is
an upper bound $q>0$ for $\kappa$ in a variational inequality. Now the question
arises, whether there is also an upper bound for $\kappa$ in an approximate
inequality. The next proposition does not answer this specific question,
but it shows that the maximal rate which can be obtained with the approach
of approximate inequalities as described in this paper is bounded by $\delta^q$.

\begin{proposition}\label{th:app_kappa_bound}
Let the assumptions of Proposition~\ref{th:kappa_bound} be satisfied, but let $u^\dagger$
satisfy an approximate inequality in the sense of Definition~\ref{def:appvari}
with $d(r)>0$ for all $r\geq 0$.
Then, with the notation of Theorem~\ref{th:appvari_rates},
$$\delta^q=\cO(\Phi^{-1}(\delta))\quad\text{as $\delta\to 0$}$$
holds.
\end{proposition}

\begin{proof}
Assume that the assertion is not true, i.e.
\begin{equation}\label{eq:misc04}
\frac{d(\Phi^{-1}(\delta))}{\delta^q}\to 0\quad\text{as $\delta\to 0$}
\end{equation}
holds.
As in the proof of Propostion~\ref{th:kappa_bound}, but starting with the inequality
$$-\xi(\tilde{u}-u^\dagger)\leq\beta_1\cB_\xi(\tilde{u},u^\dagger)
+\beta_2r^\gamma\cS(F(\tilde{u}),F(u^\dagger))^\kappa+d(r)$$
for $\tilde{u}\in M_{\bar{\alpha}}(\varrho\bar{\alpha})$ and $r\geq 0$ instead of
(\ref{eq:vari}), for $t\in(0,t_0]$ and $r\geq 0$ we get
$$-\xi(u)\leq\beta_1\underbrace{\left(\tfrac{\Omega(u^\dagger+tu)-\Omega(u^\dagger)}{t}-\xi(u)\right)}
_{\stackrel{t\to+0}{\longrightarrow}0}
+\beta_2r^\gamma\underbrace{\left(\tfrac{\cS(F(u^\dagger+tu),F(u^\dagger))^q}{t}\right)
^{\frac{\kappa}{q}}}
_{\stackrel{t\to+0}{\longrightarrow}L_\cS^{\kappa/q}}t^{\frac{\kappa}{q}-1}
+\tfrac{d(r)}{t}.$$
Now we choose $r(t):=\Phi^{-1}(t^{\frac{1}{q}})$, i.e. we have $t=\Phi(r(t))^q$.
On the one hand this (together with (\ref{eq:misc04})) implies
$$\frac{d(r(t))}{t}=\frac{d(\Phi^{-1}(t^{\frac{1}{q}}))}{t}\to 0\quad\text{as $t\to +0$}$$
and on the other hand this implies
$$r(t)^\gamma t^{\frac{\kappa}{q}-1}=\frac{r(t)^\gamma\Phi(r(t))^\kappa}{t}
=\frac{d(r(t))}{t}\to 0\quad\text{as $t\to +0$}.$$
So all terms of the above inequality tend to zero as $t\to +0$ and thus
$\xi(u)\geq 0$ holds, which is a contradiction to the assumption $\xi(u)<0$.
\end{proof}

The role of $\gamma$ in an approximate inequality is not completely clear at the moment.
If we assume that a distance function $d$ has a majorant $\bar{d}$ of the form
$\bar{d}(r)=ar^{-b\gamma}$ with $a>0$ and $b>0$ then the auxiliary functions in
Theorem~\ref{th:appvari_rates} become
\begin{equation}
\Psi(r)=a_1r^{\frac{-b\gamma(p-\kappa)-\gamma p}{\kappa}}\quad\text{and}\quad
\Phi(r)=a_2r^{\frac{-b\gamma-\gamma}{\kappa}}
\end{equation}
with constants $a_1,a_2>0$ and thus the theorem provides the convergence rate
\begin{equation}\label{eq:no_gamma_rate}
\cB_\xi(u_{\alpha(\delta)}^\delta,u^\dagger)=\cO\bigl(\delta^{\frac{b}{b+1}\kappa}\bigr),
\end{equation}
which is independent of $\gamma$.
This example shows that at least in some cases the constant $\gamma$ plays no role.
In the proof of the next proposition, however, we will see that distance functions with
majorants $\bar{d}(r)=ar^{-b\gamma}$ may occur. Furthermore, sufficient conditions for the occurrence of
some $\gamma>0$ in that context can also be found in Lemma~\ref{th:banach_rates_lemma}.

If an $\Omega$-minimizing solution satisfies a variational
inequality, then Theorem~\ref{th:vari_rates} gives us the corresponding
convergence rate. Now an interesting question is whether in this
case also an approximate inequality with higher $\kappa$ is
satisfied and, if so, does Theorem~\ref{th:appvari_rates} provide
the some rates as Theorem~\ref{th:vari_rates}? The next proposition
answers this question.

\begin{proposition}
Let $u^\dagger$ be an $\Omega$-minimizing solution which satisfies a
variational inequality in the sense of Definition~\ref{def:vari}
with $0<\kappa<p$ and let $\mu\in(\kappa,p)$ be such that $u^\dagger$
does not satisfy a variational inequality with $\kappa$ replaced by
$\mu$. Then $u^\dagger$ satisfies an approximate inequality in the
sense of Definition~\ref{def:appvari} with $\kappa$ replaced by
$\mu$ and the rates obtained from the variational inequality with
$\kappa$ and from the approximate inequality with $\mu$ coincide.
\end{proposition}

\begin{proof}
Let $\varrho$, $\bar{\alpha}$, $\beta_1$, and $\beta_2$ be the
constants from the variational inequality satisfied by $u^\dagger$;
let $\gamma>0$ be arbitrary. For all $u\in
M_{\bar{\alpha}}(\varrho\bar{\alpha})$ and all $r>0$ then
\begin{align*}
&-\xi(u-u^\dagger)
-\beta_1\cB_\xi(u,u^\dagger)-\beta_2r^\gamma\cS(F(u),F(u^\dagger))^\mu\\
&\qquad\qquad\qquad\qquad
\leq\beta_2\cS(F(u),F(u^\dagger))^\kappa-\beta_2r^\gamma\cS(F(u),F(u^\dagger))^\mu
\end{align*}
follows and
\begin{equation}
ab-\varepsilon a^{p_1}\leq\frac{b^{p_2}}{(\varepsilon p_1)^{p_2/p_1}p_2}
\end{equation}
for $a,b\geq 0$, $\varepsilon>0$, $p_1,p_2>1$ and $\frac{1}{p_1}+\frac{1}{p_2}=1$ with
\begin{align*}
&a:=\bigl(\beta_2r^\gamma\cS(F(u),F(u^\dagger))^\mu\bigr)^{\frac{\kappa}{\mu}},\quad
b:=\beta_2^{\frac{\mu-\kappa}{\mu}}r^{-\frac{\gamma\kappa}{\mu}},\\
&\varepsilon:=1,\quad p_1:=\frac{\mu}{\kappa},\quad p_2:=\frac{\mu}{\mu-\kappa}
\end{align*}
implies
\begin{align*}
d(r)
&\leq\max_{M_{\bar{\alpha}}(\varrho\bar{\alpha})}
\bigl(-\xi(u-u^\dagger)
-\beta_1\cB_\xi(u,u^\dagger)-\beta_2r^\gamma\cS(F(u),F(u^\dagger))^\mu\bigr)\\
&\leq\bigl(\tfrac{\mu}{\kappa}\bigr)^{\frac{\kappa}{\mu-\kappa}}\tfrac{\mu}{\mu-\kappa}
\beta_2r^{-\frac{\gamma\kappa}{\mu-\kappa}},
\end{align*}
i.e. $u^\dagger$ satisfies an approximate inequality and the
corresponding distance function $d$ has a majorant $\bar{d}$ of the
form $\bar{d}(r)=ar^{-b\gamma}$ with $a>0$ and
$b=\frac{\kappa}{\mu-\kappa}$.
\par Equation (\ref{eq:no_gamma_rate}) with $\kappa$ replaced by $\mu$ thus gives
\begin{equation}
\cB_\xi(u_{\alpha(\delta)}^\delta,u^\dagger)=\cO(\delta^{\frac{b}{b+1}\mu})
=\cO(\delta^{\frac{\kappa}{\mu}\mu})=\cO(\delta^\kappa)
\end{equation}
for the parameter choice
$\alpha(\delta)=c\delta^{p-\frac{b}{b+1}\mu}=c\delta^{p-\frac{\kappa}{\mu}\mu}
=c\delta^{p-\kappa}$
with a constant $c>0$.
This is exactly the convergence rate which is stated by Theorem~\ref{th:vari_rates}.
\end{proof}

\section{Source conditions and variational inequalities}\label{sec:source}

An important question which remains to be answered is the interplay
of approximate source conditions and approximate inequalities. Note
that we discussed the relationships between classical source
conditions and variational inequalities in the last paragraph of
Section~\ref{sec:notation} (see also~\cite{Hof09a}).

At first we want to show that the concept of approximate variational
inequalities described in this paper is a generalization of the
concept of approximate source conditions in Banach spaces as
introduced in \cite{HeiHof09}. So in this section our focus is on the situation of
Example~\ref{ex:banach} and we let $U$ and $V$ be reflexive
Banach spaces with $\tau_U$ and $\tau_V$ describing the corresponding
weak topologies. We set $\cS(v_1,v_2):=\Vert v_1-v_2\Vert_V$ for
$v_1,v_2\in V$, i.e. we are concerned with the Tikhonov functional
\begin{equation} \label{eq:tiknorm}
\cT_\alpha^\delta(u)=\Vert F(u)-v^\delta\Vert^p_V+\alpha\Omega(u)
\end{equation}
with $\Vert v^\delta-v^0\Vert_V\leq\delta$. In Example~\ref{ex:banach}
we mentioned that item (iv) of Assumption~\ref{as:basic} is satisfied with $s=1$. We
moreover assume that $F$, $D(F)$ and $\Omega$ are chosen such that the other items of Assumption~\ref{as:basic} are
satisfied, too.

For the remaining part of this section let $u^\dagger\in D_B$ be an $\Omega$-minimizing
solution. Because pre-compact subsets of reflexive Banach spaces are bounded
for all $\alpha$ there is a constant $K_\alpha>0$, such that
\begin{equation}
\Vert u-u^\dagger\Vert_U \leq K_\alpha\quad\text{for all $u\in
M_\alpha(\varrho\alpha)$}
\end{equation}
holds.

We make the following additional assumptions.

\begin{assumption}\label{as:banach}
It holds:
\begin{itemize}
\item[(i)]
$D(F)$ is starlike with respect to $u^\dagger$, i.e. for all $u\in D(F)$ there is a
$t_0>0$, such that $u^\dagger+t(u-u^\dagger)\in D(F)$ holds for all $t\in[0,t_0]$.
\item[(ii)]
There is a bounded linear operator $F^\prime(u^\dagger):U\rightarrow V$,
such that
$$\left\Vert\frac{F(u^\dagger+t(u-u^\dagger))-F(u^\dagger)}{t}
-F^\prime(u^\dagger)(u-u^\dagger)\right\Vert_V\to 0\quad\text{as
$t\to+0$}$$ holds for all $u\in D$.
\end{itemize}
\end{assumption}

The convexity of $\Omega$ and Assumption~\ref{as:banach} (i) imply that $D$ is then also  starlike
with respect to $u^\dagger$. In the sequel we denote by
$F^\prime(u^\dagger)^\ast:V^\ast\rightarrow U^\ast$ the adjoined operator of
$F^\prime(u^\dagger)$, where $U^\ast$ and $V^\ast$ are the dual spaces of $U$ and $V$
with respect to the norm topologies.
The handling of weakly continuous linear functionals becomes much simpler
by the fact that a linear functional on a Banach space is weakly continuous if and only if
it is continuous with respect to the norm topology.

We now define what we understand under source conditions.

\begin{definition}\label{df:source}
The $\Omega$-minimizing solution $u^\dagger$ \emph{satisfies a source
condition} if there exists an element
$\xi\in\partial\Omega(u^\dagger)$ with
$\xi\in\cR(F^\prime(u^\dagger)^\ast)$. The $\Omega$-minimizing
solution $u^\dagger$ \emph{satisfies an approximate source
condition} if there exists an element
$\xi\in\partial\Omega(u^\dagger)$ with
$\xi\in\overline{\cR(F^\prime(u^\dagger)^\ast)}$ and we define the
corresponding \emph{distance function}
$\tilde{d}:[0,\infty)\rightarrow[0,\infty)$ by
$$\tilde{d}(r):=\min\{\Vert\xi-F^\prime(u^\dagger)^\ast\eta\Vert_{U^\ast}:
\eta\in V^*,\,\Vert\eta\Vert_{V^*}\leq r\}.$$
\end{definition}

As mentioned in \cite{HeiHof09} the distance function $\tilde{d}$ is
well-defined, non-negative, finite and monotonically decreasing. If
$u^\dagger$ satisfies a source condition, then it obviously also
satisfies an approximate source condition and there is an $r_0\geq
0$ with $\tilde{d}(r)=0$ for all $r\geq r_0$. If $u^\dagger$
satisfies an approximate source condition with
$\xi\in\overline{\cR(F^\prime(u^\dagger)^\ast)}\setminus\cR(F^\prime(u^\dagger)^\ast)$
then $\tilde{d}(r)>0$ holds for all $r\geq 0$ and $\tilde{d}$ is
strictly monotonically decreasing.

The following definition was used in~\cite{HeiHof09} and
\cite{Hof09a}.

\begin{definition}
Let $c_1,c_2\geq 0$. The operator $F$ is said to be \emph{nonlinear of degree $(c_1,c_2)$}
with respect to $\Omega$, $u^\dagger$ and $\xi\in\partial\Omega(u^\dagger)$ if there
exist constants $\varrho$ fulfilling (\ref{eq:rhodef}), $\bar{\alpha}>0$, and $K>0$,
such that
$$\Vert F(u)-F(u^\dagger)-F^\prime(u^\dagger)(u-u^\dagger)\Vert_V
\leq K\Vert F(u)-F(u^\dagger)\Vert_V^{c_1}\cB_\xi(u,u^\dagger)^{c_2}$$
holds for all $u\in M_{\bar{\alpha}}(\varrho\bar{\alpha})$.
\end{definition}

The following lemma is an adaption of results in \cite{HeiHof09}.

\begin{lemma}\label{th:banach_rates_lemma}
Let the $\Omega$-minimizing solution $u^\dagger$ satisfy an
approximate source condition and let $F$ be nonlinear of degree
$(c_1,c_2)$ with respect to $\Omega$, $u^\dagger$, and $\xi$ with
$c_1\in(0,1-c_2]$ and $c_2\in[0,1)$. Further let
$\kappa:=\frac{c_1}{1-c_2}$ and $r_0>0$. Then there exist constants
$\beta_1\in[0,1)$, $\beta_2\geq 0$, and $\gamma>0$, such that
$$-\langle\xi,u-u^\dagger\rangle_{U^\ast,U}\leq\beta_1\cB_\xi(u,u^\dagger)+\beta_2r^\gamma
\Vert F(u)-F(u^\dagger)\Vert_V^\kappa+K_{\bar{\alpha}}\tilde{d}(r)$$
holds for all $u\in M_{\bar{\alpha}}(\varrho\bar{\alpha})$ and all $r\geq r_0$.
It holds $\beta_1=c_2$ and $\gamma=\frac{1}{1-c_2}$.
\end{lemma}

\begin{proof}
For $r\geq 0$ let $\eta_r\in V^\ast$ be an element for which the minimum in the definition
of $\tilde{d}(r)$ is attained. Then for $u\in M_{\bar{\alpha}}(\varrho\bar{\alpha})$ we have
\begin{align*}
\lefteqn{-\langle\xi,u-u^\dagger\rangle_{U^\ast,U}}\\
&\quad\leq\bigl\vert\langle F^\prime(u^\dagger)^\ast\eta_r+\xi-F^\prime(u^\dagger)^\ast\eta_r,
u-u^\dagger\rangle_{U^\ast,U}\bigr\vert\\
&\quad=\bigl\vert\langle\eta_r,F^\prime(u^\dagger)(u-u^\dagger)\rangle_{V^\ast,V}
+\langle\xi-F^\prime(u^\dagger)^\ast\eta_r,u-u^\dagger\rangle_{U^\ast,U}\bigr\vert\\
&\quad\leq\underbrace{\Vert\eta_r\Vert_{V^\ast}}_{\leq r}\Vert F^\prime(u^\dagger)(u-u^\dagger)\Vert_V
+\underbrace{\Vert\xi-F^\prime(u^\dagger)^\ast\eta_r\Vert_{U^\ast}}_{=\tilde{d}(r)}
\underbrace{\Vert u-u^\dagger\Vert_U}_{\leq K_{\bar{\alpha}}}\\
&\quad\leq r\Vert F(u)-F(u^\dagger)-F^\prime(u^\dagger)(u-u^\dagger)+F(u^\dagger)-F(u)\Vert_V
+K_{\bar{\alpha}}\tilde{d}(r)\\
&\quad\leq Kr\Vert F(u)-F(u^\dagger)\Vert_V^{c_1}\cB_\xi(u,u^\dagger)^{c_2}
+r\Vert F(u)-F(u^\dagger)\Vert_V+K_{\bar{\alpha}}\tilde{d}(r).
\end{align*}
Now we have to distinguish between two cases:
\begin{itemize}
\item
Case $c_2=0$. We get
\begin{align*}
\lefteqn{-\langle\xi,u-u^\dagger\rangle_{U^\ast,U}}\\
&\quad\leq Kr\Vert F(u)-F(u^\dagger)\Vert_V^{c_1}
+r\Vert F(u)-F(u^\dagger)\Vert_V+K_{\bar{\alpha}}\tilde{d}(r)\\
&\quad=\bigl(Kr+r\Vert F(u)-F(u^\dagger)\Vert_V^{1-c_1}\bigr)\Vert F(u)-F(u^\dagger)\Vert_V^{c_1}
+K_{\bar{\alpha}}\tilde{d}(r)\\
&\quad\leq\bigl(K+(\varrho\bar{\alpha})^{\frac{1-c_1}{p}}\bigr)r\Vert F(u)-F(u^\dagger)\Vert_V^{c_1}
+K_{\bar{\alpha}}\tilde{d}(r).
\end{align*}
\item
Case $c_2\in(0,1)$. We apply the inequality
$$ab\leq\frac{a^{p_1}}{p_1}+\frac{b^{p_2}}{p_2}\quad\text{for }\quad a,b\geq 0,\quad
\frac{1}{p_1}+\frac{1}{p_2}=1,\quad p_1,p_2>1$$
with
$$a:=\cB_\xi(u,u^\dagger)^{c_2},\;b:=Kr\Vert F(u)-F(u^\dagger)\Vert_V^{c_1},\;
p_1:=\frac{1}{c_2},\;p_2:=\frac{1}{1-c_2}$$
and get
\begin{align*}
\lefteqn{-\langle\xi,u-u^\dagger\rangle_{U^\ast,U}}\\
&\quad\leq Kr\Vert F(u)-F(u^\dagger)\Vert_V^{c_1}\cB_\xi(u,u^\dagger)^{c_2}
+r\Vert F(u)-F(u^\dagger)\Vert_V+K_{\bar{\alpha}}\tilde{d}(r)\\
&\quad\leq c_2\cB_\xi(u,u^\dagger)+(1-c_2)K^{\frac{1}{1-c_2}}r^{\frac{1}{1-c_2}}
\Vert F(u)-F(u^\dagger)\Vert_V^{\frac{c_1}{1-c_2}}\\
&\quad\quad\,+r\Vert F(u)-F(u^\dagger)\Vert_V+K_{\bar{\alpha}}\tilde{d}(r)\\
&\quad=c_2\cB_\xi(u,u^\dagger)+K_{\bar{\alpha}}\tilde{d}(r)\\
&\quad\quad\,+\Bigl((1-c_2)K^{\frac{1}{1-c_2}}r^{\frac{1}{1-c_2}}
+r\Vert F(u)-F(u^\dagger)\Vert_V^{\frac{1-c_1-c_2}{1-c_2}}\Bigr)
\Vert F(u)-F(u^\dagger)\Vert_V^{\frac{c_1}{1-c_2}}\\
&\quad\leq c_2\cB^\cR_\xi(u,u^\dagger)+K_{\bar{\alpha}}\tilde{d}(r)\\
&\quad\quad\,+\Bigl((1-c_2)K^{\frac{1}{1-c_2}}
+\varrho_{\bar{\alpha}}^{\frac{1-c_1-c_2}{p(1-c_2)}}r_0^{\frac{-c_2}{1-c_2}}\Bigr)
r^{\frac{1}{1-c_2}}\Vert F(u)-F(u^\dagger)\Vert_V^{\frac{c_1}{1-c_2}}.\qedhere
\end{align*}
\end{itemize}
\end{proof}

\begin{theorem}\label{th:d_tilde_d}
Let the $\Omega$-minimizing solution $u^\dagger$ satisfy an
approximate source condition and let $F$ be nonlinear of degree
$(c_1,c_2)$ with respect to $\Omega$, $u^\dagger$ and $\xi$ with
$c_1\in(0,1-c_2]$, $c_2\in[0,1)$, and $\frac{c_1}{1-c_2}<p$. Then
$u^\dagger$ satisfies an approximate inequality in the sense of
Definition~\ref{def:appvari} with $0<\kappa=\frac{c_1}{1-c_2}<p$ and
$$d(r)\leq K_{\bar{\alpha}}\tilde{d}(r)\quad\text{for all $r\geq r_0>0$}$$
with $r_0$ from Lemma~\ref{th:banach_rates_lemma} holds.
\end{theorem}

\begin{proof}
The assertion is a direct consequence of Lemma~\ref{th:banach_rates_lemma}.
\end{proof}

If $u^\dagger$ satisfies a source condition then by Theorem~\ref{th:d_tilde_d}
$u^\dagger$ also satisfies a variational inequality
and Theorem~\ref{th:vari_rates} and \cite[Theorem 3.3]{HeiHof09}
provide the same convergence rate. In analogy we have: If
$u^\dagger$ satisfies an approximate source condition then $u^\dagger$
also satisfies an approximate inequality and the rates obtained in
Theorem~\ref{th:appvari_rates} are not worse than the rates in
\cite[Theorem 4.3]{HeiHof09}.

\begin{remark}
In \cite{HofKalPoeSch07} and \cite{HeiHof09} it has been shown that
in the case $c_1=0$ and $c_2=1$ a source condition
$\xi=F^\prime(u^\dagger)^\ast\eta$ with $K\Vert\eta\Vert_{V^*}<1$ implies
a variational inequality with $\kappa=1$. The converse result that a
variational inequality with $\kappa=1$ implies the source condition
is true if $F$ and $\Omega$ are G\^ateaux differentiable in
$u^\dagger$ (see \cite{SchGraGroHalLen09}). However, the authors think
that convergence results in the case  $c_1=0$ and $c_2=1$ are missing
when $u^\dagger$ only satisfies an approximate source condition in the sense of
Definition~\ref{df:source} with $\tilde d(r)>0$ for all $r \ge 0.$
\end{remark}

Now that we know about a basic relationship between approximate source conditions
and variational inequalities we conclude this section by repeating from \cite{Hof09a} the interplay of
classical H\"older type source conditions and variational inequalities in Hilbert spaces.
So let $U$ and $V$ be Hilbert spaces and let $F=A$ be a bounded linear
operator.
Taking the standard Tikhonov functional
$$\cT_\alpha^\delta(u)=\Vert Au-v^\delta\Vert^2_V+\alpha\Vert u\Vert^2_U$$
the subdifferential of $\Omega=\Vert\mybullet\Vert^2_U$ at $u\in U$
is the singleton $\{\langle\mybullet,2u\rangle_U\}$
(where $\langle\mybullet,\mybullet\rangle_U$ denotes the inner product), i.e. we set
$\xi=\langle\mybullet,2u^\dagger\rangle_U$, and the corresponding Bregman distance
is $\cB_\xi(\mybullet,u^\dagger)=\Vert\mybullet-u^\dagger\Vert_U^2$. To legitimize the
extended concept of variational inequalities for $\kappa\neq 1$ in
\cite{Hof09a} the following is stated:
%Taking the standard Tikhonov functional
%$$\cT_\alpha^\delta(u)=\Vert Au-v^\delta\Vert^2_V+\alpha\Vert u\Vert^2_U$$
%the subdifferential of $\Omega=\Vert\mybullet\Vert^2_U$ at $u\in U$
%is the singleton $\{\langle\mybullet,2u\rangle_U\}$, and we set
%$\xi=\langle\mybullet,2u^\dagger\rangle_U$. To legitimize the
%extended concept of variational inequalities for $\kappa\neq 1$ in
%\cite{Hof09a} the following is stated:

If $u^\dagger$ satisfies a source condition of type
$u^\dagger\in\cR((A^\ast A)^{\frac{\mu}{2}})$ with $\mu\in(0,1)$
then $u^\dagger$ satisfies a variational inequality
\begin{equation} \label{eq:vilin}
\langle u^\dagger-u,2u^\dagger\rangle_U\leq
\beta_1\Vert u-u^\dagger\Vert_U^2+\beta_2\Vert A(u-u^\dagger)\Vert_V^\kappa
%\langle\xi,u^\dagger-u\rangle_{U^\ast,U} \le
%\beta_1\cB_\xi(u,u^\dagger)+\beta_2
%\|A(u-u^\dagger)\|_V^\kappa
\end{equation}
with
$\kappa=\frac{2\mu}{1+\mu}$. For $\mu=1$ this holds too, as we saw
in the Banach space setting above. Because of
Proposition~\ref{th:kappa_bound} such a relationship cannot hold for
$\mu>1$.
In \cite[Proposition 5.7]{Hof09a} also the following converse result is formulated:
If $u^\dagger$ satifies a variational inequality (\ref{eq:vilin}) with exponent $\kappa$ then it
satisfies a source condition of type $u^\dagger\in\cR((A^\ast A)^{\frac{\mu}{2}})$
for all $\mu\in(0,\frac{\kappa}{2-\kappa})$.

\section{Conclusions and open questions} \label{sec:conclude}

The following diagram should help to understand the cross-connections between the different approaches
for obtaining convergence rates. In this context, $\Rightarrow$ stands for an implication and
$\rightarrow$ stands for ``as good or better as''. This, however, is only a very rough characterization
of the interplay which the reader can find in detail in the corresponding theorems,
propositions and remarks.

$$
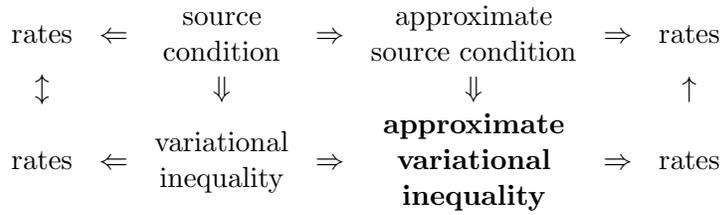
\begin{array}{c}
\begin{array}{ccccccc}
\text{rates}&\Leftarrow&\begin{matrix}\text{source}\\
\text{condition}\end{matrix}&\Rightarrow&\begin{matrix}\text{approximate}\\
\text{source condition}\end{matrix}&\Rightarrow&\text{rates}\\
\updownarrow&&\Downarrow&&\Downarrow&&\uparrow\\
\text{rates}&\Leftarrow&\begin{matrix}\text{variational}\\
\text{inequality}\end{matrix}&\Rightarrow&\begin{matrix}\text{\bf approximate}\\
\text{\bf variational}\\\text{\bf inequality}\end{matrix}&\Rightarrow&\text{rates}\\
\end{array}\\
\end{array}
$$

\bigskip

As we have seen from Proposition~\ref{th:kappa_bound}, Remark~\ref{rem:rem1},
and Proposition~\ref{th:app_kappa_bound} for the Banach space setting
when (\ref{eq:normex}) and (\ref{eq:tiknorm}) are under consideration
the proven convergence rates of Section~\ref{sec:vari} and Section~\ref{sec:appvari}
are because of the occurring limitation $\kappa \le 1$
by construction not faster than $\cB_\xi(u_{\alpha(\delta)}^\delta,u^\dagger)=\cO(\delta)$ as
$\delta\to 0$. Therefore with the technique of variational inequalities (\ref{eq:vinenorm})
and also with the corresponding approximate inequalities we are
captured in the \emph{low rate world}. A \emph{higher rate world} for that Banach space setting was structured, for example, by
the recent papers \cite{Hein09,Neub09}, where under higher source conditions, for $p>1$,  and provided that the space $V$ is smooth enough rates up to
$\cB_\xi(u_{\alpha(\delta)}^\delta,u^\dagger)=\cO(\delta^{4/3})$ can be proven.

In our low rate world the rates are additionally limited by the inequality $\kappa<p$.
Up to now the literature considered preferably the case $p > 1$, where this inequality
gives no restriction. In the case $0<p \le 1$, however, for which our approach also applies, this gives a
serious restriction. One can interpret the condition $\kappa<p$ then as follows: The exponent
$0<p < 1$ seems to be a  \emph{qualification} of the chosen method (similar to the qualification of
linear regularization methods, see \cite{EHN96}) which itself defines an upper bound for
convergence rates. If the \emph{smoothness} of the solution $u^\dagger$ grows further, i.e.~$p<\kappa \le 1$, then
the convergence rate does not follow. Note that the boundary situation $0<\kappa=p \le 1$ shows the
so-called \emph{exact penalization} effect studied in \cite{BurOsh04} for $p=1$, where the rate
$\cB_\xi(u_{\alpha_{\mathrm{fix}}}^\delta,u^\dagger)=\cO(\delta)$  was proven under the source condition $\xi=F^\prime(u^\dagger)^*\eta$ whenever the regularization parameter $\alpha_{\mathrm{fix}}>0$ was chosen fixed but small enough.
>From the proof of Lemma~\ref{th:rates_lemma} yielding the estimate (\ref{eq:pgleichkappa}) we immediately obtain the corresponding rate
 $\cB_\xi(u_{\alpha_{\mathrm{fix}}}^\delta,u^\dagger)=\cO(\delta^p)$  whenever a
 variational inequality is satisfied with exponent $0<p=\kappa \le 1$  and the regularization parameter $\alpha_{\mathrm{fix}}>0$ is fixed and small enough.
However, it is an open problem to answer the question whether the  rates
$\cO(\delta^{\min\{\kappa,p\}})$ for $0<p < 1$  can be improved or not.

An advantage of our new approach for the low rate world is the fact that the items (ii) and (iii) of
Proposition~\ref{prop:levelsets} tell us that $\{M_\alpha(\varrho\alpha):\alpha>0\}$ in some sense is a family of
neighbourhoods of solutions $u^\dagger$ to the equation $F(u)=v^0$. We recall that if a variational inequality
holds on a level set $M_{\bar{\alpha}}(\varrho\bar{\alpha})$ then it holds on each
level set $M_{\alpha}(\varrho\alpha)$ with $0<\alpha<\bar{\alpha}$. Hence, satisfying
a variational inequality means that there exists an arbitrarily small neighbourhood
of $u^\dagger$ such that a variational inequality holds on this neighbourhood.
Or, in other words, convergence rates depend only on the behaviour of the three
functionals $\xi(\mybullet-u^\dagger)$, $\cB_\xi(\mybullet,u^\dagger)$ and
$\cS(F(\mybullet),F(u^\dagger))$ in an arbitrarily small neighborhood of the set of
solutions. Looking at the problem from such \emph{functional} point of view
this suggests the conjecture that some kind of variational inequality like tools may exist,
which is able to integrate higher source conditions and would lead us to the higher rate world.
For example, we see that $\cS$ and $F$ themselves are  not important, only their
combination $\cS(F(\mybullet),F(u^\dagger))$ is of interest. Hence one could ask in this context
how the mentioned functionals
reflect the combination of source conditions and structure of nonlinearity in case
of higher smoothness. This should be forthcoming work.

\bibliography{geissler_hofmann}

\end{document}